\documentclass[final,3p,times,twocolumn]{elsarticle}
\usepackage{amsmath}

\usepackage{amssymb}
\usepackage{amsmath,amsfonts,bm,mathtools}
\usepackage{algorithmicx}
\usepackage{algpseudocode}
\usepackage{cases}
\usepackage{multicol}
\usepackage{color}

\newtheorem{remark}{Remark}[section]
\newtheorem{example}{Example}[section]

\journal{Neurocomputing}

\begin{document}

\begin{frontmatter}



\title{Discrete state transition algorithm for unconstrained \\ integer optimization problems}


\author{Xiaojun Zhou$^{\dag,\ddag}$, David Yang Gao$^{\ddag}$, Chunhua Yang$^{\dag}$, Weihua Gui$^{\dag}$}

\address{$^{\dag}$School of Information Science and Engineering, Central South University, Changsha 410083, China\\
$^{\ddag}$School of Science, Information Technology and Engineering, Federation University Australia, Victoria 3353, Australia.

}

\begin{abstract}
A recently new intelligent optimization algorithm called discrete state transition algorithm is considered in this study, for solving unconstrained integer optimization problems. Firstly, some key elements for discrete state transition algorithm are summarized to guide its well development. Several intelligent operators are designed for local exploitation and global exploration. Then, a dynamic adjustment strategy ``risk and restoration in probability" is proposed to capture global solutions with high probability. Finally, numerical experiments are carried out to test the performance of the proposed algorithm compared with other heuristics, and they show that the similar intelligent operators can be applied to ranging from traveling salesman problem, boolean integer programming, to discrete value selection problem, which indicates the adaptability and flexibility of the proposed intelligent elements.
\end{abstract}


\end{frontmatter}


\section{Introduction}
\label{intro}
In this paper, we consider the following unconstrained integer optimization problem
\begin{eqnarray}\label{eqn1}
\min f(x),
\end{eqnarray}
where, $x = (x_1, \cdots, x_n) \in \mathbb{Z}^n$. \\
\indent Generally speaking, the above optimization problem is NP-hard, which can not be solved in polynomial time. A direct method is to adopt the so called ``divide-and-conquer" strategy, which separates the optimization problem into several subproblems and then solve these subproblems step by step. Branch and bound (B\&B), branch and cut (B\&C), and branch and price (B\&P) belong to this kind; however, these methods are essentially in exponential time.  An indirect method is to relax the optimization problem by loosening its integrality constraints to continuity and then solve the continuous relaxation problem or its Lagrangian dual problem, including LP-based relaxation, SDP-based relaxation, Lagrangian relaxation, etc. Nevertheless, when rounding off the relaxation solution, they may cause some infeasibility or can only get approximate solutions, and when using Lagrangian dual, there may exist duality gap between the primal and the dual problem \cite{dowman,geoffrion,junger,li}. \\
\indent On the other hand, some stochastic algorithms, such as genetic algorithm (GA) \cite{ahmed,yokota}, simulated annealing (SA) \cite{kirkpatrick,zhang}, ant colony optimization (ACO) \cite{dorigo,schluter}, are also widely used for integer optimization problems, which aim to obtain ``good solutions" in reasonable time. In terms of the concepts of state and state transition, a new heuristic search algorithm called state transition algorithm (STA) has been proposed recently, which exhibits excellent global search ability in continuous function optimization \cite{xzhou2011a,xzhou2011b,xzhou2012,xzhou2013,xzhou2014}. In \cite{yang}, three intelligent operators (geometrical operators) named swap, shift and symmetry have been designed for discrete STA to solve the traveling salesman problem (TSP), and it shows that the discrete STA outperforms its competitors with respect to both time complexity and search ability. In \cite{xzhou2015}, a discrete state transition algorithm is successfully applied to the optimal design of water distribution networks.
 To better develop discrete STA for medium-size or large-size discrete optimization problems, in the study, we firstly build the framework of discrete state transition algorithm and propose five key elements for discrete STA, of which, the representation of a decision variable, the local and global operators and the dynamic adjustment strategy are mainly studied. Four geometrical operators named swap, shift, symmetry and substitute are designed, which are intelligent due to their adaptability and flexibility in various types of integer optimization.
The mixed strategies of ``greedy criterion" and ``risk and restoration in probability" are proposed, in which, ``greedy criterion" and ``restoration in probability" are used to guarantee a good convergence performance, and  ``risk a bad solution in probability" aims to escape from local optimality. Some applications ranging from traveling salesman problem, boolean integer programming, to discrete value selection problem are studied. Experimental results have demonstrated the effectiveness and efficiency of the proposed method.

The main contribution and novelty of this paper is three-fold, which can be summarized as follows:
1) a systematic formulation of discrete state transition algorithm is firstly proposed, including the
state space representation and five key elements;
2) a dynamic adjustment strategy called ``risk and restoration in probability" is designed to improve the ability to
escape from local optima;
3) the proposed algorithm is successfully integrated with several classical integer optimization problems.

\section{The framework of discrete state transition algorithm}
If a solution to a specific optimization problem is described as a state, then the transformation to update the solution becomes a state transition. Without loss of generality, the unified form of generation of solution in discrete state transition algorithm can be described as
\begin{eqnarray}
\left \{ \begin{array}{ll}
x_{k+1}= A_{k}(x_{k}) \bigoplus B_{k}(u_{k})\\
y_{k+1}= f(x_{k+1})
\end{array} \right.,
\end{eqnarray}
where, $x_{k} \in \mathbb{Z}^{n}$ stands for a current state, corresponding to a solution of a specific optimization problem; $u_{k}$ is a function of $x_{k}$ and historical states; $A_{k}(\cdot)$, $B_{k}(\cdot)$ are transformation operators, which are usually state transition matrixes; $\bigoplus$ is a operation, which is admissible to operate on two states; $f$ is the fitness function.\\
\indent As an intelligent optimization algorithm, discrete state transition algorithm have the following five key elements:\\
\indent (1) Representation of a solution. In discrete STA, we choose a special representation, that is, the permutation of the set $\{1,2,\cdots,n\}$, which can be easily manipulated by some intelligent operators. The reason that we call the operators ``intelligent" is due to their geometrical property (swap, shift, symmetry and substitute), and an intelligent operator has the same geometrical function for different types of problems. A big advantage of such a representation and operators is that, after each state transformation, the newly created state is always feasible, avoiding the trouble into rounding off a continuous solution into an integral one.\\
\indent (2) Sampling in a candidate set. When a transformation operator is exerted on a current state, the next state is not deterministic, that is to say, there are possibly different choices for the next state. It is not difficult to imagine that all possible choices will constitute a candidate set, or a ``neighborhood". Then we execute several times of transformation, called search enforcement (\textit{SE}) degree, on current state, to sample in the ``neighborhood". Sampling is a very important factor in state transition algorithm, which can characterize the search space and avoid enumeration.\\
\indent (3) Local exploitation and global exploration. In continuous optimization, it is quite significant to design good local and global operators. The local exploitation can guarantee high precision of a solution and convergent performance of a algorithm, and the global exploration can avoid getting trapped into local minima or prevent premature convergence.
In discrete optimization, it is extremely difficult to define a ``good" local optimal solution due to its dependence on a problem's structure, which leads to the same difficulty in the definition of local exploitation and global exploration. Anyway, in discrete state transition algorithm, we define a little change to current solution by a transformation as local exploitation, while a big change to current solution by a transformation as global exploration.\\
\indent (4) Self learning and regular communication. State transition algorithm behaves in two styles, one is individual-based, the other is population-based, which is certainly a extended version. The individual-based state transition algorithm focuses on self learning, in other words, it focuses on designing operators and dynamic adjustment (details given in the following). Undoubtedly, communication among different states is a promising strategy for state transition algorithm, as indicated in \cite{xzhou2012}. Through communication, states can share information and cooperate with each other. However, how to communicate and when to communicate are key issues. In continuous state transition algorithm, intermittent exchange strategy was proposed, which means that states communicate with each other at a certain frequency in a regular way.\\
\indent (5) Dynamic adjustment. It is a potentially useful strategy for state transition algorithm. In the iterative process of searching, the fitness value can decrease sharply in the early stage, but it stagnates in the late stage, due to the static environment. As a result, some perturbation should be added to activate the environment. In fact, dynamic adjustment can be understood and implemented in various ways. For example, the alternative use of different local and global operators is a dynamic adjustment to some extent. Then, we can change the search enforcement degree, vary the fitness function, reduce the dimension, etc. Of course, ``risk a bad solution in probability" is another dynamic adjustment, which is widely used in simulated annealing (SA). In SA, the Metropolis criterion \cite{metroplis} is used to accept a bad solution:
$p = \mathrm{exp}(\frac{-\triangle E}{k_B T})$,
where, $\triangle E = f(\bm x_{k+1}) - f(\bm x_{k})$, $k_B$ is the Boltzmann probability  factor, $T$ is the temperature to regulate the process of annealing. In the early stage, temperature is high, and it has big probability to accept a bad solution, while in the late stage, temperature is low, and it has very small probability to accept a bad solution, which is the key point to guarantee the convergence. We can see that the Metropolis criterion has the ability to escape from local optimality, but on the other hand, it will miss some ``good solutions" as well.

In discrete STA, a novel strategy, named ``risk and restoration in probability", is proposed. Details can be found in the following individual-based STA.

\subsection{Individual-based discrete STA}
\indent In this part, we focus on the individual-based discrete STA, and the main process of discrete STA is shown in the pseudocode as follows
\begin{algorithmic}[1]
\Repeat
    \State {[Best,fBest] $\gets$ swap(*,Best,fBest)}
    \State {[Best,fBest] $\gets$ shift(*,Best,fBest)}
    \State {[Best,fBest] $\gets$ symmetry(*,Best,fBest)}
    \State {[Best,fBest] $\gets$ substitute(*,Best,fBest)}
    \If{ fBest $<$ fBest$^{*}$} \Comment{greedy criterion}
    \State {Best$^{*}$ $\gets$ Best}
    \State {fBest$^{*}$ $\gets$ fBest}
    \EndIf
    \If{$rand < p_1$} \Comment{restoration in probability}
    \State {Best $\gets$ Best$^{*}$}
    \State {fBest $\gets$ fBest$^{*}$}
    \EndIf
\Until{the specified termination criteria are met}
\end{algorithmic}
\quad \\
\indent To be specific, swap function in above pseudocode is given as follows for example
\begin{algorithmic}[1]
\State{State $\gets$ op\_swap(*,Best)}
\State{[newBest,fnewBest] $\gets$ fitness(*,State)}
\If{fnewBest $<$ fBest}\Comment{greedy criterion}
    \State{Best $\gets$ newBest}
    \State{fBest $\gets$ fnewBest}
\Else
    \If{$rand < p_2$}\Comment{risk in probability}
    \State{Best $\gets$ newBest}
    \State{fBest $\gets$ fnewBest}
    \EndIf
\EndIf
\end{algorithmic}

\begin{figure}[!htbp]
  \centering
  \includegraphics[width=8cm,height = 10cm]{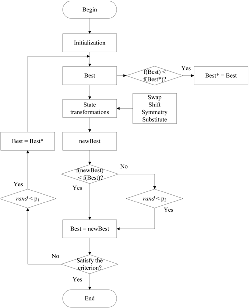}
  \caption{The flowchart of the individual-based discrete STA}\label{flowchart}
\end{figure}

\indent From the pseudocode, we can find that, on the whole, ``greedy criterion" is adopted to keep the incumbent ``\textit{Best}$^{*}$"; while partially in the inner process, a bad solution ``\textit{Best}" is accepted in each state transformation at a  probability $p_2$, and in the same time, the ``\textit{Best}$^{*}$" is restored in the outer process at another probability $p_1$. The ``risk a bad solution in probability" strategy aims to escape from local optimality, while the ``greedy criterion" and ``restoring the incumbent best solution in probability" are to guarantee a good convergence performance. The flowchart of the individual-based dynamic STA is illustrated in Fig.\ref{flowchart}, and we can find that
the incumbent best ``\textit{Best}$^{*}$" is kept in an external archive.

\subsection{Population-based discrete STA}
As indicated in \cite{xzhou2012}, a population-based approach can improve the performance of STA to a large extent.
The crossover operation is a typical way for individual states to communicate with each other. Let $x_1$ and $x_2$ be individual components of old states, $y_1$ and $y_2$ are components of new states. In this paper, a simple crossover operator
is inherited as follows
\begin{equation}
\left \{ \begin{array}{ll}
y_1= \delta x_1 +(1-\delta) x_2,\\
y_2= (1-\delta) x_1 + \delta x_2,
\end{array} \right.
\end{equation}
where, $\delta$ is a random variable, which obeys the 0-1 distribution.

\begin{figure}[!htbp]
  \centering
  \includegraphics[width=5cm,height = 5cm]{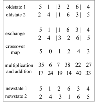}
  \caption{Tie-breaking crossover in population-based discrete STA}\label{crossover}
\end{figure}

The above crossover operation can be utilized directly for many cases; however, for traveling salesman problem, a repairing procedure is necessary to generate a
feasible solution. In this study, we introduce the tie-breaking crossover \cite{poon}, and the procedure can be found in Fig.\ref{crossover},
in which, the crossover map is a random ordering of the integers
$0,1,\cdots,n-1$; the intermediate states are multiplied by $n$ and added to the corresponding numbers in the crossover map; after a sort procedure, the new states are created at last.
The pseudocode of the kernel of population-based discrete STA can be outlined in the following
\begin{algorithmic}[1]
\State {[State,fState] $\gets$ initiation(*)}
\State {BState $\gets$ State, BfState $\gets$ fState}
\Repeat
    \State {[State,fState,BState,BfState] $\gets$ self\_learning()}
    \If{mod(iter,CF)==0}
        \State {[State,fState] $\gets$ communication(*)}
    \EndIf
    \If{fState $<$ BfState}
        \State {BState $\gets$ State, BfState $\gets$ fState}
    \EndIf
    \If{$rand$ $<$ $p_1$}
        \State {State $\gets$ BState, fState $\gets$ fBState}
    \EndIf
    \State {[Best$^{*}$,fBest$^{*}$] $\gets$ fitness(*,BState)}
\Until{the specified termination criteria are met}
\end{algorithmic}
where, \textit{self\_learning} means each state will be performed on several state transformations (swap, shift, etc), which is similar to the
individual-based discrete STA; \textit{communication} corresponds to the
crossover operation, and \textit{CF} is the communication frequency. By the way, $``*"$ stands for some omitted parameters such as search enforcement (\textit{SE}),
number of states (\textit{SN}).
\begin{remark}
It should be noted that there are two main differences between individual-based discrete STA and individual-based metaheuristics as well as
population-based discrete STA and population-based metaheuristics.
The first one is how to generate candidate solutions. In discrete STA, four special state transformation operators are designed, called swap, shift, symmetry and substitute, respectively. Those transformation operators have different geometrical properties and they can make sure that, after each state transformation, the newly created state is always feasible, which avoids the trouble into rounding off a continuous solution into an integral one as other metaheuristics.
Another important difference in generating solutions is the sampling mechanism used in discrete STA.
In theory, the number of a complete candidate solution set for a state transformation operator can be very large;
however, to avoid enumeration, we execute several times of state transformation (called search enforcement (\textit{SE}) degree), to sample in the solution set.
The second one is how to accept a new solution. In discrete STA, a novel accepting criterion called ``risk and restoration in probability" is proposed. In the inner process, a relatively bad solution is accepted in probability, while in the outer process, the historical best solution is restored in probability, which aims to jump out of local minima as well as to maintain good convergence results.
\end{remark}

\section{Theoretical analysis of the discrete STA}
In this section, we analyze the convergence performance, global search ability, and time complexity of the discrete state transition algorithm.

Similarly to the continuous case, we give the definition of local and global minima for unconstrained integer optimization as follows
\begin{subequations}
\begin{eqnarray}\label{eq3a}
f(x^{*})  \leq f(x), \forall\; x \in \mathbb{Z}^n,
\end{eqnarray}
\begin{eqnarray}\label{eq3b}
f(x^{*}) \leq  f(x), \forall\; | x^{*} - x | < 1, x^{*}, x \in \mathbb{Z}^n
\end{eqnarray}
\end{subequations}
If (\ref{eq3a}) is satisfied, we say that $x^{*}$ is a global minimizer, and if (\ref{eq3b}) is satisfied, the $x^{*}$ is called a local minimizer.

\textbf{Theorem 1} The sequence generated by discrete STA can converge to a local minimizer.\\
\indent \textit{Proof}.
Let us suppose the maximum number of iterations (denoted by \textit{M}) is big enough, considering that the ``greedy criterion" is used to keep the incumbent best $Best^{*}_k$, then we have
$
f(Best^{*}_{k+1}) \leq f(Best^{*}_k)
$, that is to say, the sequence $\{f(Best^{*}_k)\}$ is a monotonically decreasing sequence,
and
there must exist a number $N < M$, when $k > N$, updating of the incumbent best will no longer happen, i.e., $f(Best^{*}_k) = f(Best^{*}_N), \forall ~ k > N$, where $Best^{*}_N$ is the solution in the \textit{N}th iteration.
On the other hand, by the definition of a local minimizer in (\ref{eq3b}), we can find that every integral solution is a local minimizer. Let $Best^{*}_N$ denote as the local minimum solution $x^{*}$, and then we have $f(Best^{*}_N) - f(x^{*}) = 0$. \qed

To show the discrete STA converges in probability to the set of globally minimum states,
let $\mathbb{Z}_0$ denote the set of states in $\mathbb{Z}^n$ at which $f(\cdot)$ attains its global minimum value,
\begin{eqnarray*}
\mathbb{Z}_0 &=& \{x \in \mathbb{Z}^n| f(x) - f(x^{*}) < \varepsilon\}, \; \forall  \varepsilon > 0 \\
\mathbb{Z}_1 &=& \mathbb{Z}^n \backslash \mathbb{Z}_0
\end{eqnarray*}
and assume that any state in $\mathbb{Z}^n$ is reachable from any other state in $\mathbb{Z}$.

\textbf{Theorem 2} The sequence generated by discrete STA can converge to the global minimum in probability.\\
\indent \textit{Proof}.
It is obvious to find that
the random process $Best = (Best^{*}_k: k \geq 0)$ produced by the discrete STA is a discrete time
Markov chain. The one-step transition probability matrix at step $k$ is
\begin{eqnarray*}
P(Best^{*}_{k+1} \in \mathbb{Z}_0 |Best^{*}_{k} \in \mathbb{Z}_0) &=& 1 \\
P(Best^{*}_{k+1} \in \mathbb{Z}_1 |Best^{*}_{k} \in \mathbb{Z}_0) &=& 0 \\
P(Best^{*}_{k+1} \in \mathbb{Z}_0 |Best^{*}_{k} \in \mathbb{Z}_1) &\geq& c \\
P(Best^{*}_{k+1} \in \mathbb{Z}_1 |Best^{*}_{k} \in \mathbb{Z}_1) &\leq& 1 - c
\end{eqnarray*}
where, $c$ is the lower bound of the transition probability from $\mathbb{Z}_1$ to $\mathbb{Z}_0$.
Due to the assumption of reachability, we can find that the discrete Markov chain is irreducible and $c \in (0,1)$.
By using the similar methodology of Markov ergodic convergence theorem in \cite{hajek1988},
we have
\begin{eqnarray*}
\lim_{k \rightarrow \infty} P(Best^{*}_{k} \in \mathbb{Z}_0) = 1.
\end{eqnarray*}
\begin{remark}
The global search ability depends on the assumption of reachability to a large extent.
To meet the assumption required by the analysis, two fundamental elements exist in discrete STA.
The first one is related to the global operators, which have the functionality of
bringing a big change to current solution. Another is the `risk in probability' strategy
in dynamic discrete STA since a relative bad solution is accepted in probability.
Anyhow, the theoretical convergence result requires that the the number of iterations approaches to infinity.
In practice, some additional strategies and techniques need to be
added to improve its practical search ability.
\end{remark}

With respect to the time complexity of discrete STA, it should be noted that
the discrete STA aims to obtain a satisfactory solution in a reasonable amount of time.
In the above pseudocodes as described, it can be found that in the outer loop, there are \textit{M} iterations, while in the inner loop, there exist four times of \textit{SE} transformations, that is to say, the time complexity of the proposed discrete STA is $\mathcal{O}(M \cdot SE)$.

Next, some applications are given to describe the details, from the traveling salesman problem, boolean integer programming, to discrete value selection problem.

\section{Application for traveling salesman problem}
Suppose $\mathcal{N} = \{1, \cdots, n\}$ is the set of cities, the traveling salesman problem (TSP) can be described as: given a set of \textit{n} cities and the distance $d_{ij}$ for each pair of cities $i$ and $j$, find a roundtrip of minimal total length visiting each city exactly once. Typically, the traveling salesman problem is usually modeled as the following two representations \cite{smith}.\\
\indent (LP-TSP):
\begin{eqnarray}
~~~ & \min\limits_{x_{ij}} ~~\sum\limits_{i=1}^{n} \sum\limits_{j=1}^{n}x_{ij}d_{ij} \nonumber \\
s.t. & \sum\limits_{i=1}^{n} x_{ij} = 1, \forall j \in \mathcal{N} \nonumber \\
     & \sum\limits_{j=1}^{n} x_{ij} = 1, \forall i \in \mathcal{N} \nonumber \\
     & \sum\limits_{i \in S}  \sum\limits_{j \in \bar{S}} x_{ij} \geq 1, \nonumber \\
     & \forall S \subset \mathcal{N}, \bar{S} \subset \mathcal{N} \setminus S \nonumber \\
     & x_{ij} \in \{0,1\}, \forall i, j \in  \mathcal{N},
\end{eqnarray}
where, the decision variable $x_{ij}$ is defined by
\begin{subnumcases}
{x_{ij}=}
1, & if city \textit{i} is followed by city \textit{j}\\
0, & otherwise
\end{subnumcases}
\indent (QP-TSP):
\begin{eqnarray}
~~~ & \min\limits_{x_{ij}} ~~\sum\limits_{i=1}^{n} \sum\limits_{k=1}^{n} \sum\limits_{j=1}^{n}x_{ij}d_{ik}(x_{k(j+1)} + x_{k(j-1)}) \nonumber \\
s.t. & \sum\limits_{i=1}^{n} x_{ij} = 1, \forall j \in \mathcal{N} \nonumber \\
     & \sum\limits_{j=1}^{n} x_{ij} = 1, \forall i \in \mathcal{N} \nonumber \\
     & x_{ij} \in \{0,1\}, \forall i, j \in  \mathcal{N},
\end{eqnarray}
here, the decision variable $x_{ij}$ is defined by
\begin{subnumcases}
{x_{ij}=}
1, & if city \textit{i} is in the \textit{j}th position \\
0, & otherwise
\end{subnumcases}
\indent We can find that the model of linear programming (LP) based TSP is different from the model of quadratic programming (QP) based TSP in two aspects. One is the definition of the decision variable $x_{ij}$, the other is that (LP-TSP) has one more constraint than (QP-TSP). Taking the difficulty of dealing with constraints into consideration, in discrete STA, we use another simple representation which can be easily manipulated by intelligent operators.
\subsection{State transformation operators for TSP}
\begin{figure}[!htbp]
  \centering
  \includegraphics[width=6.5cm]{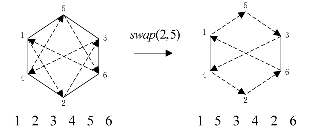}
  \caption{Illustration of swap transformation for TSP}\label{swaptsp}
\end{figure}

As for a \textit{n}-city traveling salesman problem, a permutation of $\{1,2, \cdots, n\}$ is used to represent a solution to the problem in discrete STA. Based on the representation, three special transformation operators are proposed to illustrate local and global search.\\
\indent (1) swap transformation
\begin{eqnarray}
\bm x_{k+1}= A^{swap}_{k}(m_a) \bm x_{k},
\end{eqnarray}
where, $A^{swap}_{k} \in \mathbb{R}^{n \times n}$ is called swap transformation matrix, $m_a$ is a constant integer called swap factor to control the maximum number of positions to be exchanged, while the positions are random. If $m_a = 2$, we call the swap operator local exploitation, and if $m_a \geq 3$, the swap operator is regarded as global exploration in this case. Fig.\ref{swaptsp} gives the function of the swap transformation graphically when $m_a = 2$. In this case, the state transition process is as follows
\begin{eqnarray*}
\begin{scriptsize}
\left( {\begin{array}{*{20}{c}}
   1\\
   5\\
   3\\
   4\\
   2\\
   6\\
\end{array}} \right) =
\left( {\begin{array}{*{20}{c}}
   1 & 0 & 0 & 0 & 0 & 0\\
   0 & 0 & 0 & 0 & 1 & 0\\
   0 & 0 & 1 & 0 & 0 & 0\\
   0 & 0 & 0 & 1 & 0 & 0\\
   0 & 1 & 0 & 0 & 0 & 0\\
   0 & 0 & 0 & 0 & 0 & 1\\
\end{array}} \right)
\times
\left( {\begin{array}{*{20}{c}}
   1\\
   2\\
   3\\
   4\\
   5\\
   6\\
\end{array}} \right)
\end{scriptsize}
\end{eqnarray*}

\begin{figure}[!htbp]
  \centering
  \includegraphics[width=6.5cm]{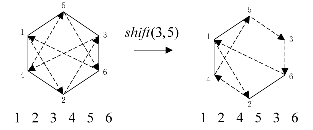}
  \caption{Illustration of shift transformation for TSP}\label{shifttsp}
\end{figure}

(2) shift transformation
\begin{eqnarray}
\bm x_{k+1}= A^{shift}_{k}(m_b) \bm x_{k},
\end{eqnarray}
where, $A^{shift}_{k} \in \mathbb{R}^{n \times n}$ is called shift transformation matrix, $m_b$ is a constant integer called shift factor to control the maximum length of consecutive positions to be shifted. By the way, the selected position to be shifted after and positions to be shifted are chosen randomly. Similarly, shift transformation is called local exploitation and global exploration when $m_b = 1$ and $m_b \geq2 $  respectively. To make it more clearly, if $m_b = 1$, we set position 3 to be shifted after position 5, as described in Fig.\ref{shifttsp}. In this case, we have
\begin{eqnarray*}
\begin{scriptsize}
\left( {\begin{array}{*{20}{c}}
   1\\
   2\\
   4\\
   5\\
   3\\
   6\\
\end{array}} \right) =
\left( {\begin{array}{*{20}{c}}
   1 & 0 & 0 & 0 & 0 & 0\\
   0 & 1 & 0 & 0 & 0 & 0\\
   0 & 0 & 0 & 1 & 0 & 0\\
   0 & 0 & 0 & 0 & 1 & 0\\
   0 & 0 & 1 & 0 & 0 & 0\\
   0 & 0 & 0 & 0 & 0 & 1\\
\end{array}} \right)
\times
\left( {\begin{array}{*{20}{c}}
   1\\
   2\\
   3\\
   4\\
   5\\
   6\\
\end{array}} \right)
\end{scriptsize}
\end{eqnarray*}

\begin{figure}[!htbp]
  \centering
  \includegraphics[width=6.5cm]{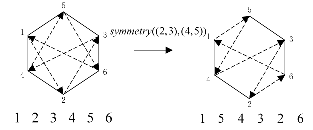}
  \caption{Illustration of symmetry transformation for TSP}\label{symtsp}
\end{figure}

(3) symmetry transformation
\begin{eqnarray}
\bm x_{k+1}= A^{sym}_{k}(m_c) \bm x_{k},
\end{eqnarray}
where, $A^{sym}_{k} \in \mathbb{R}^{n \times n}$ is called symmetry transformation matrix, $m_c$ is a constant integer called symmetry factor to control the maximum length of subsequent positions as center. By the way, the component before the subsequent positions and consecutive positions to be symmetrized are both created randomly. Considering that the symmetry transformation can make big change to current solution, it is intrinsically called global exploration. For instance, if $m_c = 0$, let choose component 3, then the subsequent position or the center is \{$\emptyset$\}, the consecutive positions are $\{4,5\}$, and the function of symmetry transformation is given in Fig.\ref{symtsp}. Then, we have
\begin{eqnarray*}
\begin{scriptsize}
\left( {\begin{array}{*{20}{c}}
   1\\
   5\\
   4\\
   3\\
   2\\
   6\\
\end{array}} \right) =
\left( {\begin{array}{*{20}{c}}
   1 & 0 & 0 & 0 & 0 & 0\\
   0 & 0 & 0 & 0 & 1 & 0\\
   0 & 0 & 0 & 1 & 0 & 0\\
   0 & 0 & 1 & 0 & 0 & 0\\
   0 & 1 & 0 & 0 & 0 & 0\\
   0 & 0 & 0 & 0 & 0 & 1\\
\end{array}} \right)
\times
\left( {\begin{array}{*{20}{c}}
   1\\
   2\\
   3\\
   4\\
   5\\
   6\\
\end{array}} \right)
\end{scriptsize}
\end{eqnarray*}
\onecolumn\begin{multicols}{2}
\begin{remark}
The transformation matrix comes from the identity matrix by
random elementary transformation. Taking the swap operator for example, the swap transformation matrix
can be computed by the following pseudocode
\begin{verbatim}
function y = swap_matrix(n,ma)
y  = eye(n);
R  = randperm(n);
T  = R(1:ma);
S  = T(randperm(ma));
y(T,:) = y(S,:);
\end{verbatim}
The above swap matrix has the effect to exchange any two random positions for a solution.
\end{remark}
\subsection{TSP instances for test}
To evaluate the performance of the individual-based discrete STA (DSTAI) and the population-based discrete STA (DSTAII) as well as the previously proposed discrete STA (DSTA0) in \cite{yang} for traveling salesman problem, some medium-size TSP instances are used for test.
We set $m_a = 2, m_b = 1, m_c = 0, m_d = 1, p_1 = 0.1, p_2 = 0.1$, \textit{SN} = 10, \textit{SE} = 20, \textit{CF} = 1
and the maximum number of iterations at 1500.

The compared SA is a combination of a recently published one from \cite{geng} and the other from \cite{seshadri}, in which, there are three mutations, namely, vertex insert (VI), block insert (BI) and block reverse (BR),
and the ratios of VI, BI and BR in the proposed algorithm were designed for $10\%, 1\%$ and $89\%$, respectively.
The initial temperature for the SA is $2000$, and the cooling rate is $0.97$.
Since SA (simulated annealing) is individual based algorithm, the maximum number of iterations of SA is set at 90000 for fair comparison.

Programs are run independently for 20 trails for each algorithms in MATLAB R2010b (version of 7.11.0.584) on Intel(R) Core(TM) i3-2310M CPU @2.10GHz under Window 7 environment. Some statistics as well as the ``error" are computed for comparison, where, ``error" is defined by
\begin{eqnarray}
error = \frac{best - optimum}{optimum} \times 100\%, \nonumber
\end{eqnarray}
here, $best$ is the best result achieved by discrete STA or SA, $optimum$ is the incumbent best result in TSPLIB \cite{reinelt}. Experimental results for these TSP instances are given in Table \ref{tspinstances}.
\end{multicols}
\begin{table*}[!htbp]
\centering
\caption{Experimental results for TSP instances}
\label{tspinstances}
\begin{tabular}{{ccccccc}}
\hline\noalign{\smallskip}
\textit{instance} & \textit{optimum}    & \textit{algorithm} &\textit{best} & \textit{mean} & \textit{s.t.}  & \textit{error}\\
\noalign{\smallskip}\hline\noalign{\smallskip}
kroA100.tsp& 21282    & SA        & 2.1729e4   &2.2635e4 & 778.7240  & 2.10\% \\
(\textit{n}=100)  &   & DSTA0     & 2.1853e4   &2.3213e4 & 906.1100  & 2.69\% \\
           &          & DSTAI     & 2.1782e4   &2.2835e4 & 715.8493  & 2.35\% \\
           &          & DSTAII    & 2.1294e4   &2.1767e4 & 221.6416  & 0.0583\% \\
\hline
kroB100.tsp& 22141    & SA        & 2.3032e4   &2.3657e4 & 445.7826  & 4.03\%\\
(\textit{n}=100)   &  & DSTA0     & 2.3230e4   &2.3794e4 & 517.0476  & 4.92\%\\
           &          & DSTAI     & 2.3012e4   &2.3734e4 & 507.3792  & 3.93\% \\
           &          & DSTAII    & 2.2345e4   &2.2880e4 & 302.1363  & 0.9229\% \\
\hline
kroC100.tsp& 20749    & SA        & 2.1417e4   &2.2223e4 & 522.2034  & 3.22\%\\
(\textit{n}=100)  &   & DSTA0     & 2.1275e4   &2.2877e4 & 709.8698  & 2.53\%\\
           &          & DSTAI     & 2.1038e4   &2.1891e4 & 536.8809  & 1.40\% \\
           &          & DSTAII    & 2.0907e4   &2.1378e4 & 246.3382  & 0.76\% \\
\hline
kroD100.tsp& 21294    & SA        & 2.1896e4   &2.2911e4 & 483.0088  & 2.83\%\\
(\textit{n}=100)   &  & DSTA0     & 2.1945e4   &2.3043e4 & 565.7970  & 3.06\%\\
           &          & DSTAI     & 2.1867e4   &2.2665e4 & 592.5252  & 2.69\% \\
           &          & DSTAII    & 2.1380e4   &2.1991e4 & 315.3214  & 0.40\% \\
\hline
kroE100.tsp& 22068    & SA        & 2.2523e4   &2.3125e4 & 389.4191  & 2.06\%\\
(\textit{n}=100)  &   & DSTA0     & 2.2692e4   &2.3738e4 & 450.8241  & 2.83\% \\
           &          & DSTAI     & 2.2419e4   &2.3371e4 & 678.6940  & 1.59\% \\
           &          & DSTAII    & 2.2311e4   &2.2637e4 & 166.8205  & 1.10\% \\
\noalign{\smallskip}\hline
\end{tabular}
\end{table*}

\begin{multicols}{2}
\indent As can be seen from Table \ref{tspinstances}, the DSTAI has almost the same performance as the SA,
while the DSTAII gets the best results, the biggest error of which is almost no more than 1\%,
which testifies the effectiveness of the proposed operators and strategies, namely, the ``risk and restoration in probability" can be comparable with the
Metropolis criterion in SA.

\begin{remark}
It should be noted that the results of SA in this paper are different from those of SA in \cite{geng} for the same instances.
The reason is that we only use the mutation operators proposed in \cite{geng}, combing with a standard SA in \cite{seshadri},
removing the two-stage adaptive local search strategy.
\end{remark}

\section{Application for boolean integer programming}
In boolean integer programming (BIP), a solution is comprised of a series of boolean values ($\mathcal{I} = \{0,1\}$ or $\mathcal{I} = \{-1,1\}$). Swap, shift and symmetry operators can also be applied to internal transformation (operators aiming to change the internal components of a sequence), and another operator called substitute is designed for external transformation (operator aiming to bring alien components into a sequence). It should be noted that $\mathcal{I} = \{0,1\}$ is the same to $\mathcal{I} = \{-1,1\}$ under such a circumstance, although there exists a linear transformation relationship between them in other studies.
\subsection{State transformation operators for BIP}
\begin{figure*}[!htbp]
  \centering
  \includegraphics[width=12cm,height = 1.2cm]{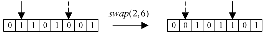}
  \caption{Illustration of swap transformation for BIP}
  \label{bswap}
\end{figure*}

\begin{figure*}[!htbp]
  \centering
  \includegraphics[width=12cm,height = 1.5cm]{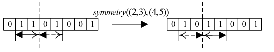}
  \caption{Illustration of symmetry transformation for BIP}
  \label{bsymmetry}
\end{figure*}

As we mentioned previously, the same intelligent operator has the same geometrical property for different applications. It is not difficult to imagine that the swap, shift and symmetry operators for boolean integer programming have the same formulation as that of traveling salesman problem. Let $\mathcal{I} = \{0,1\}$, the illustrations of internal transformation are given from Fig.\ref{bswap} to Fig.\ref{bsymmetry}.

\begin{figure*}[!htbp]
  \centering
  \includegraphics[width=12cm,height = 1.2cm]{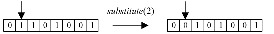}
  \caption{illustration of substitute transformation for BIP}
  \label{bsubstitute}
\end{figure*}

Next, let introduce the external transformation.\\
\indent (4) substitute transformation
\begin{eqnarray}
\bm x_{k+1}= A^{sub}_{k}(m_d) \bm x_{k},
\end{eqnarray}
where, $A^{sub}_{k} \in \mathbb{R}^{n \times n}$ is called substitute transformation matrix, $m_d$ is a constant integer called substitute factor to control the maximum number of positions to be substituted. By the way, the positions are randomly created. If $m_d = 1$, we call the substitute operator local exploitation, and if $m_d \geq 2$, the substitute operator is regarded as global exploration in this case. Fig.\ref{bsubstitute} gives the function of the substitute transformation vividly when $m_d = 1$.
\end{multicols}
\begin{figure*}[!htbp]
  \centering
  \includegraphics[width=12cm,height = 1.2cm]{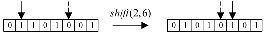}
  \caption{Illustration of shift transformation for BIP}
\end{figure*}

\begin{multicols}{2}
\subsection{Maxcut instances for test}
Let $G = (V,E)$ be an undirected graph with edge weight $w_{ij}$ on $n+1 = |V|$ vertices and $m = |E|$ edges, for each edge $(i,j) \in E$, the maximum cut problem (Maxcut) is to find a subset $S$ of the vertex set $V$ such that the total weight of the edges between $S$ and its complementary subset $\bar{S} = V \backslash S$ is as large as possible.\\
\indent LP based Maxcut model \cite{krishnan}:\\
\indent Considering a variable $y_{ij}$ for each edge $(i,j) \in E$, and assuming $y_{ij}$ to be 1 if $(i,j)$ is in the cut, and 0 otherwise, the Maxcut can be modeled as the following linear programming (LP) optimization problem:
\begin{eqnarray}\label{maxcut1}
    \max ~W(\bm y) & = & \sum_{i=1}^{n+1}\sum_{i<j,(i,j) \in E} w_{ij} y_{ij} \nonumber \\
         \mathrm{s.t.} ~ & & \bm y \mathrm{~is~the~incidence~vector~of~a~cut}, \nonumber \\
\end{eqnarray}
Here the incidence vector $\bm y = \{y_{ij}\} \in \mathbb{R}^m$, where the $m$ is the number of edges in the graph.\\
\indent Let CUT($G$) denote the convex hull of the incidence vectors of cuts in $G$. Since maximizing a linear function over a set of points equals to maximizing it over the convex hull of this set of points, we can rewrite (\ref{maxcut1}) to the following
\begin{eqnarray}
    \max ~W(\bm y) & = & \bm c^T \bm y \nonumber \\
         \mathrm{s.t.} & & \bm y \in \mathrm{CUT}(G).
\end{eqnarray}
where $c = \{w_{ij}\} \in \mathbb{R}^m$.\\
\indent QP based Maxcut model \cite{wang}:\\
\indent For a bipartition $(S,\bar{S})$, with $y_i = 1$ if $i \in S$, and $y_i = -1$ otherwise, the Maxcut can also be formulated as the following integer optimization problem:
\begin{eqnarray}\label{maxcut2}
    \max ~W(\bm y) & = &\frac{1}{4} \sum_{i=1}^{n+1}\sum_{j=1}^{n+1} w_{ij}(1 - y_i y_j) \nonumber \\
    \mathrm{s.t.} ~ & & \bm y \in\{-1,1\}^{n+1}.
\end{eqnarray}
\indent Without loss of generality, if we fix the value of the last variable at 1, then the problem (\ref{maxcut2}) is equivalent to the integer quadratic programming problem
\begin{eqnarray}\label{maxcut3}
\min\Big\{P(\bm x) = \frac{1}{2} \bm x^T Q \bm x  -  \bm x^T \bm c: \bm x \in \{-1,1\}^{n}\Big\},
\end{eqnarray}
where, $Q = \{Q_{ij}\}$ is a symmetric matrix with $Q_{ij} = w_{ij}(i,j = 1,2,\cdots, n)$, and $\bm c = -(w_{1(n+1)}, \cdots, w_{n(n+1)})^T$.  It is not difficult to find that a optimal solution $\bm x^{*}$ to problem (\ref{maxcut3}) corresponds to a optimal solution $(\bm x^{*},1)$ of original problem (\ref{maxcut2}).\\

\end{multicols}

\begin{table*}[!htbp]
\centering
\caption{Experimental results for Maxcut instances}
\label{maxcutinstances}
\begin{tabular}{{ccccccc}}
\hline\noalign{\smallskip}
\textit{instance} & \textit{optimum}    & \textit{algorithm} &\textit{best} & \textit{mean} & \textit{s.t.}  & \textit{error}\\
\hline
kroA100    & 5897392  & GA        &   5897392    &5.8622e6  & 4.5060e4  & 0 \\
(\textit{n}=100)&     & DSTAI     &   5897392    &5.8879e6  & 2.9123e4  & 0 \\
           &          & DSTAII    &   5897392    &5897392   & 2.8666e-9 & 0 \\
\hline
kroB100    & 5763047  & GA         &   5763047    &5.7444e6 & 2.4007e4  & 0\\
(\textit{n}=100)&     & DSTAI      &   5763047    &5.7529e6 & 2.0726e4  & 0\\
           &          & DSTAII     &   5763047    &5763047  & 1.9110e-9 & 0\\
\hline
kroC100    & 5890760  & GA         &   5890760    &5.8678e6 & 3.3939e4  & 0\\
(\textit{n}=100)&     & DSTAI      &   5890760    &5.8706e6 & 3.1574e4  & 0 \\
           &          & DSTAII     &   5890760    &5890760  & 0         & 0 \\
\hline
kroD100    & 5463250  & GA         &  5463250     &5.4387e6 & 3.3594e4  & 0 \\
(\textit{n}=100)&     & DSTAI      &  5463250     &5.4410e6 & 3.4850e4  & 0\\
           &          & DSTAII     &  5463250     &5463250  & 2.8666e-9 & 0\\
 \hline
kroE100    & 5986591  & GA         &  5986591     &5.9372e6 & 5.9586e4  & 0 \\
(\textit{n}=100)&     & DSTAI      &  5986591     &5.9585e6 & 4.9985e4  & 0\\
           &          & DSTAII     &  5986591     &5986591  & 9.5552e-10& 0\\
\noalign{\smallskip}\hline
\end{tabular}
\end{table*}
\begin{multicols}{2}
\indent Considering that the QP based Maxcut model is easier to manipulate for intelligent operators, it is adopted in discrete STA for simple representation.
We use a real coded integer genetic algorithm \cite{deep} for comparison, in which, laplace crossover (the location parameter $a = 0$ and the scaling parameter $b = 0.35$) and power mutation (the index of mutation $p = 4$) were used, and a truncation procedure was applied to make sure the integrity of a solution. More specifically, the crossover probability $p_c = 0.8$ and the mutation probability $p_m = 0.005$.
The parameters setting for discrete STA is the same to that in TSP instances except the maximum number of iterations at $200$ and \textit{CF} $ = 20$.
The population size for integer genetic algorithm (GA) is the same to the search enforcement (\textit{SE}) in STA, while the maximum iterations
for GA is set at four times as that of discrete STA. In the same way, we define the following ``error"
\begin{eqnarray}
error = \frac{optimum - best}{optimum} \times 100\%, \nonumber
\end{eqnarray}
here, $best$ is the best result achieved by discrete STA or integer GA, $optimum$ can be found in \cite{wang}. Experimental results for Maxcut instances are given in Table \ref{maxcutinstances}.

As can be seen from Table \ref{maxcutinstances}, all of these algorithms have the ability to achieve the global minimum for all of the instances. The performance of GA and DSTAI is much the same, while DSTAII is much superior than its competitors since it can find the global minimum for every instance in each run.
\begin{remark}
When using the integer GA for the Maxcut problem, we need to reformulate the model by
a linear transformation, namely
\begin{eqnarray}
\bm x = 2 \bm z - 1, 0 \leq \bm z \leq 1, \bm z \in Integer
\end{eqnarray}
\end{remark}
\section{Application for discrete value selection}
Typically, the formulation of discrete value selection (DVS) problem is different from the model in (\ref{eqn1}), because the domain is defined as follows
\begin{eqnarray}
x_i \in \mathcal{U} = \{u_1, \cdots, u_m\}, u_j \in \mathbb{R}, j = 1, \cdots, m.
\end{eqnarray}
By introducing a linear transformation
\begin{eqnarray}
x_i = \sum_{j=1}^{K} u_j y_{ij},
\end{eqnarray}
where
\begin{eqnarray}
\sum_{j=1}^{K} y_{ij} = 1, y_{ij} \in \{0,1\},
\end{eqnarray}
then the discrete value selection can be rewritten to the equivalent constrained boolean integer programming problem \cite{yu}.\\
\indent In discrete STA, we only use the index of $u_j$ to represent a solution, for example, a solution $(1, 3, 2)$ is corresponding to $(u_1, u_3, u_2)$, which is easy to be manipulated by the intelligent operators.

\subsection{Transformation operators for discrete value selection}
\begin{figure*}[!htbp]
  \centering
  \includegraphics[width=12cm]{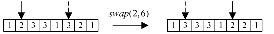}
  \caption{Illustration of swap transformation for DVS}
  \label{dvsswap}
\end{figure*}

\begin{figure*}[!htbp]
  \centering
  \includegraphics[width=12cm,height = 1.5cm]{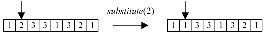}
  \caption{Illustration of substitute transformation for DVS}
  \label{dvssubstitute}
\end{figure*}

\indent The intelligent operators swap, shift, symmetry and substitute for discrete value selection are similar to that in boolean programming problem. As a result, only
illustrations of these transformations are given from Fig.\ref{dvsswap} to Fig.\ref{dvssubstitute}.
\end{multicols}

\begin{figure*}[!htbp]
  \centering
  \includegraphics[width=12cm]{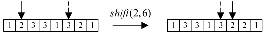}
  \caption{Illustration of shift transformation for DVS}
\end{figure*}
\begin{figure*}[!htbp]
  \centering
  \includegraphics[width=12cm,height = 2cm]{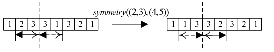}
  \caption{Illustration of symmetry transformation for DVS}
\end{figure*}
\begin{multicols}{2}
\subsection{The integer Rosenbrock function for test}
The continuous Rosenbrock function has been widely studied as a benchmark. Considering that the Rosenbrock function also has many integer local minima, in this paper, it is the first time to use it as a benchmark for integer optimization problem. The integer Rosenbrock function is defined as
\begin{eqnarray*}
f(\bm x) = \sum_{i=1}^{n-1}[100(x_{i+1} - x_i^2)^2 + (x_i - 1)^2],  \\
x_i \in \{-2,-1,0,1,2\}.
\end{eqnarray*}
It is not difficult to find that $x = (x_1, \cdots, x_n)$ for any $x_i = 0, i = 1, \cdots, n$ is a local minimum of the integer Rosenbrock function.\\
\indent The parameters setting is the same to that of Maxcut instances except the maximum number of iterations, which is specified at 10 times of the problem dimension. The population size for integer GA is the same to \textit{SE}, while the maximum iterations of GA is four times as that of discrete STA.
Experimental results for integer Rosenbrock function are given in Table \ref{irosenbrock}.
\end{multicols}

\begin{table*}[!htbp]
\centering
\caption{Experimental results for integer Rosenbrock function}
\label{irosenbrock}
\begin{tabular}{{ccccccc}}
\hline\noalign{\smallskip}
\textit{n} & \textit{optimum}    & \textit{algorithm} &\textit{best} & \textit{mean} & \textit{s.t.}  & \textit{error}\\
\hline
5       & 0   & GA        &   0        &0        & 0          & 0\\
        &     & DSTAI     &   0        &0        & 0          & 0\\
        &     & DSTAII    &   0        &0        & 0          & 0\\
\hline
10      & 0   & GA        &   0        &0        & 0          & 0\\
        &     & DSTAI     &   0        &0        & 0          & 0 \\
        &     & DSTAII    &   0        &0        & 0          & 0 \\
\hline
20      & 0  & GA         &   0        &0        & 0          & 0\\
        &    & DSTAI      &   0        &0        & 0          & 0 \\
        &    & DSTAII     &   0        &0        & 0          & 0 \\
\hline
50      & 0  & GA         &   0        &20.1000  & 61.8665    & 0\\
        &    & DSTAI      &   0        &0        & 0          & 0 \\
        &    & DSTAII     &   0        &0        & 0          & 0 \\
\hline
100     & 0  & GA         &   0        &50.2500  & 89.2966    & 0\\
        &    & DSSTA      &   0        &0        & 0          & 0 \\
        &    & DSTAII     &   0        &0        & 0          & 0 \\
\hline
200      & 0  & GA        &   0        &140.7500 & 132.0705   & 0\\
         &    & DSTAI     &   0        &0        & 0          & 0 \\
         &    & DSTAII    &   0        &0        & 0          & 0 \\
\hline
500      & 0  & GA        &   1508     &2.6594e3 & 1.0091e3   & -\\
         &    & DSTAI     &   0        &22.5000  & 46.3119    & 0 \\
         &    & DSTAII    &   0        &0        & 0          & 0 \\
\noalign{\smallskip}\hline
\end{tabular}
\end{table*}

\begin{multicols}{2}
As can be seen from Table \ref{irosenbrock}, all of these algorithms do well in the low dimensional cases. For large scale problem, discrete STA is much superior to integer GA, but the DSTAI begins to have poor performance when the dimension is larger than 500. On the other hand, the results obtained by DSTAII are much more satisfactory.

\begin{remark}
 When using the integer GA for the integer Rosenbrock problem, we need to revise the constraints to
\begin{eqnarray}
-2 \leq \bm x \leq 2, \bm x \in Integer
\end{eqnarray}
\end{remark}

\subsection{Other examples}
We continue to test some other discrete value selection examples, which can be found in \cite{wu2011}.
The maximum iterations are 100, 500 and 1000 for the following three examples respectively in discrete STA.
For integer GA, the maximum iterations is four times that of discrete STA.
For these examples, we define the following ``error"
\begin{eqnarray}
error = \frac{|optimum - best|}{|optimum|} \times 100\%, \nonumber
\end{eqnarray}
where, $best$ is the best result achieved by discrete STA or integer GA.\\
\begin{example}
\begin{eqnarray*}
&&\min f_1(x) = \frac{1}{2} \bm x^T Q \bm x + \bm c^T \bm x \nonumber \\
&&\mathrm{s.t.} \;\;\bm x \in \{0,1,2,\cdots,10\}^8
\end{eqnarray*}
\end{example}
where,
\begin{eqnarray*}
Q =
\begin{pmatrix}
     4  &  -2  &   -3    &  0     & 1   &   4   &   5   &  -2\\
    -2   &  -4   &   0   &   0   &   2   &   2    &  0   &   0\\
    -3    &  0   &   8   &  -2   &   0    &  3   &   4   &   0\\
     0   &   0    & -2    & -4   &   4   &   4   &   0   &   1\\
     1   &   2   &   0    &  4   & 100   &   2   &   0  &   -2\\
     4   &   2   &   3   &   4    &  2   & 100    &  1   &   0\\
     5   &   0  &    4   &   0    &  0    &  1   & 200   &   4\\
    -3   &   0   &   0   &   1   &  -2   &   0   &   4   &  10\\
\end{pmatrix},\\
\bm c^T =
\begin{pmatrix}
 -4  &   1  &  -8  &   3 & -100 &  -10  & -20  &   0
\end{pmatrix}.
\end{eqnarray*}
\begin{example}
\begin{eqnarray*}
&&\min f_2(x) = \bm x^T Q \bm x \nonumber \\
&&\mathrm{s.t.} \;\;\bm x \in \{0,1,2,\cdots,49\}^{10}
\end{eqnarray*}
\end{example}
where,
\begin{equation*}
Q \!=\!
\begin{pmatrix}
    -1   & -2   &  2   &  8  &  -5  &   1  &  -4   &  0 &    0  &   8\\
    -2  &   2   &  0   & -5  &   4  &  -4  &  -4   & -5  &   0  &  -5\\
     2 &    0  &   2  &  -3 &    7  &   0  &  -3  &   7  &   5   &  0\\
     8 &   -5  &  -3  &  -1  &  -3   & -1   &  7   &  1  &   7  &   2\\
    -5  &   4  &   7  &  -3   &  1  &   0  &  -4  &   2  &   4 &   -2\\
     1 &   -4 &    0  &  -1   &  0 &    1   &  9  &   5   &  2  &   0\\
    -4 &   -4  &  -3  &   7  &  -4 &    9  &   3  &   1  &   2  &   0\\
     0 &   -5  &   7 &    1  &   2 &    5  &   1  &   0  &  -3   & -2\\
     0 &    0  &   5  &   7   &  4  &   2  &   2  &  -3   &  2 &    3\\
     8  &  -5  &   0   &  2   & -2 &    0  &   0  &  -2  &   3  &   3\\
\end{pmatrix}
\end{equation*}

\begin{example}
\begin{eqnarray*}
&&\min f_3(x) = \bm x^T Q \bm x + \bm c^T \bm x \nonumber \\
&&s.t. \;\;\bm x \in \{0,1,2,\cdots,99\}^{20}
\end{eqnarray*}
\end{example}
\end{multicols}
where,
\begin{eqnarray*}
Q =
\left( \begin{array}{cccccccccccccccccccc}
    -3   &  7  &   0  &  -5  &   1   &  1  &   0   &  2 &   -1  &  -1  &  -9   &  3  &   5  &   0  &   0  &   1 &    7  &  -7   & -4  &  -6\\
     7   &  0 &   -5   &  1 &    1  &   0  &   2  &  -1  &  -1 &   -9  &   3  &   5  &   0  &   0  &   1   &  7 &   -7  &  -4  &  -6  &  -3\\
     0  &  -5  &   1   &  1  &   0  &   2  &  -1  &  -1  &  -9   &  3   &  5   &  0  &   0 &    1  &   7   & -7  &  -4  &  -6 &   -3  &   7\\
    -5 &    1 &    1  &   0  &   2  &  -1   & -1  &  -9  &   3   &  5  &   0 &    0  &   1  &   7 &   -7 &   -4  &  -6 &   -3  &   7  &   0\\
     1   &  1  &   0  &   2  &  -1 &   -1 &   -9  &   3 &    5   &  0   &  0  &   1  &   7  &  -7 &   -4  &  -6   & -3 &    7    & 0   & -5\\
     1  &   0 &    2 &   -1  &  -1  &  -9  &   3  &   5  &   0   &  0   &  1 &    7 &   -7  &  -4   & -6  &  -3 &    7 &    0  &  -5   &  1\\
     0  &   2 &   -1  &  -1  &  -9   &  3  &   5 &    0  &   0   &  1  &   7  &  -7  &  -4 &   -6 &   -3  &   7 &    0 &   -5  &   1 &    1\\
     2 &   -1  &  -1  &  -9  &   3 &    5  &   0   &  0  &   1  &   7   & -7 &   -4 &   -6  &  -3 &    7 &    0  &  -5 &    1   &  1   &  0\\
    -1  &  -1  &  -9 &    3  &   5 &    0   &  0   &  1   &  7  &  -7  &  -4  &  -6  &  -3  &   7  &   0 &   -5  &   1  &   1   &  0  &   2\\
    -1  &  -9   &  3 &    5  &   0  &   0   &  1  &   7 &   -7  &  -4   & -6 &   -3&     7   &  0 &   -5  &   1  &   1&     0  &   2  &   1\\
    -9  &   3 &    5 &    0  &   0  &   1 &    7  &  -7  &  -4  &  -6  &  -3 &    7&     0   & -5   &  1 &    1  &   0 &    2    & 1 &    2\\
     3  &   5  &   0 &    0  &   1   &  7 &   -7  &  -4  &  -6  &  -3  &  7  &   0   & -5   &  1   &  1 &    0&     2 &
     1 &    2   &  3\\
     5   &  0  &   0 &    1 &    7  &  -7  &  -4  &  -6 &   -3  &   7 &    0&    -5    & 1   &  1  &   0 &    2  &   1&     2  &   3  &   9\\
     0  &   0   &  1 &    7  &  -7  &  -4  &  -6  &  -3   &  7   &  0   & -5   &  1  &   1   &  0  &   2   &  1 &    2 &    3   &  9   &  4\\
     0  &   1  &   7 &   -7  &  -4  &  -6  &  -3   &  7  &   0  &  -5  &   1  &   1  &   0 &    2  &   1  &   2   &  3 &    9   &  4 &   -1\\
     1  &   7  &  -7 &   -4  &  -6  &  -3  &   7  &   0  &  -5    & 1   &  1 &    0   &  2   &  1 &    2 &    3 &    9&     4   & -1 &   -3\\
     7 &   -7  &  -4  &  -6  &  -3  &   7  &   0 &   -5  &   1  &   1  &   0   &  2 &    1  &   2 &    3  &   9  &   4 &   -1 &   -3  &   9\\
    -7  &  -4  &  -6 &   -3  &   7  &   0  &  -5   &  1  &   1 &    0  &   2   &  1  &  2   &  3   &  9   &  4  & -1   & -3 &    9  &   7\\
    -4  &  -6  &  -3  &   7  &   0   & -5  &   1    & 1   &  0   &  2 &    1   &  2  &   3  &   9  &   4  &  -1&    -3  &   9    & 7   & -9\\
    -6  &  -3  &   7  &   0  &  -5  &   1 &    1  &   0  &   2 &    1 &    2   &  3  &  9 &   4   & -1   & -3&     9   &  7
    & -9  &   8\\
\end{array}\right)
\end{eqnarray*}
and
\begin{eqnarray*}
\bm c^T =
\left( \begin{array}{cccccccccccccccccccc}
 -5   &  2  &  -1   & -3   &  5   &  4 &   -1   &  0  &   9  &   4   &  7   & -4   &  3 &    5 &    8   & -1  &   1  &   5    & -6  &   9
\end{array}\right)
\end{eqnarray*}

\begin{table*}[!htbp]
\centering
\caption{Experimental results for other DVS examples}
\label{otherexamples}
\begin{tabular}{{ccccccc}}
\hline\noalign{\smallskip}
instances &\textit{optimum}    & \textit{algorithm} &\textit{best} & \textit{mean} & \textit{s.t.}  & \textit{error}\\
$f_1$& -620    & GA    & -620     & -616.6750   & 12.8106        & 0\\
     &         & DSTAI & -620     & -620        & 0              & 0\\
     &         & DSTAII& -620     & -620        & 0              & 0\\
$f_2$& -70429  & GA    & -70429   & -6.4980e4   & 6.1412e3       & 0\\
     &         & DSTAI & -70429   & -6.9909e4   & 2.3264e3       & 0\\
     &         & DSTAII& -70429   & -70429      & 0              & 0\\
$f_3$& -1439658& GA    & -1407590 & -1.2397e6   & 9.9237e4       & 2.23\%\\
     &         & DSTAI & -1439658 & -1.3486e6   & 9.8920e4       & 0\\
     &         & DSTAII& -1439658 & -1.3871e6   & 6.9660e4       & 0\\
\noalign{\smallskip}\hline
\end{tabular}
\end{table*}

\begin{multicols}{2}
As can be seen from Table \ref{otherexamples}, the discrete STA outperforms integer GA for all the three cases on the whole. For the former two examples, all of these algorithms have the capability to obtain the global optimum.  However, for the last example, only discrete STA can achieve the global optimum, which indicates its superiority to its competitor, and it also shows
that population-based STA has better performance than individual-based STA.
\section{Conclusion}
In this paper, a new intelligent optimization algorithm named discrete state transition algorithm is studied for integer optimization problem. It is the first time to build the framework for
discrete state transition algorithm, and five key elements are discussed to better develop the algorithm. The representation of a feasible solution and the dynamic adjustment strategy are mainly studied. Various applications have shown the adaptability and flexibility of the designed intelligent operators and experimental results have testified the effectiveness and efficiency of the proposed algorithm and strategies.

On the other hand, it should be noted that the proposed discrete STA is not good enough, especially for the traveling salesman problem.
In the future, we will extend the proposed discrete STA to more efficient ones by using additional strategies to
improve the global search ability and reduce the computational cost.
\section*{Acknowledgements}
We would also like to thank the anonymous
reviewers for their valuable comments and suggestions that
helped improve the quality of this paper.
This work was supported by the National Science Foundation for Distinguished Young Scholars of China
(61025015), the Foundation for Innovative Research Groups of the National Natural Science
Foundation of China (61321003), the
National Natural Science Foundation of China (Grant No.
61503416) and the State Key Program of National Natural Science
of China (Grant Nos.61533020 and 61533021).




\bibliographystyle{model1a-num-names}
\bibliography{<your-bib-database>}



\end{multicols}

\end{document}